\documentclass[11pt]{article} 

\usepackage{amsmath}

\title{A note on doubles of groups}
\author {Nadia Benakli \thanks{e-mail: \sl benakli@math.columbia.edu} \and
         Oliver T. Dasbach \thanks{e-mail: \sl kasten@math.columbia.edu, 
         \sl http://www.math.uni-duesseldorf.de/home/kasten} \and
         Yair Glasner \thanks{e-mail: \sl yair@math.columbia.edu}
         \and 
         Brian Mangum \thanks{e-mail: \sl mangum@math.columbia.edu}}

\numberwithin{equation}{section}

\def\defeq{\ensuremath{\stackrel{\mathrm{def}}{=}}} 

\usepackage{amsthm, amssymb, amscd}

\usepackage[dvips]{graphicx}

\newcommand {\Z} {\mathbb{Z}}
\newcommand{\inclusion}{\hookrightarrow}
\newcommand{\arrow}{\rightarrow}
\newcommand{\subgroup}{<}
\newcommand{\infinity}{\infty}

\newtheorem{theorem}{Theorem}[section]
\newtheorem{proposition}[theorem]{Proposition}
\newtheorem{lemma}[theorem]{Lemma}
\newtheorem{corollary}[theorem]{Corollary}
\newtheorem{definition}[theorem]{Definition}

\date{}

\begin{document}

\maketitle

\begin{abstract}
Recently, Rips produced an example of a double of two free groups which
has unsolvable generalized word problem.  In this paper, we show that
Rips's example fits into a large class of doubles of groups, each member
of which contains $F_2 \times F_2$ and therefore has unsolvable
generalized word problem and is incoherent. 
\end{abstract}

\addtocounter{section}{1}

By a result of Mihailova every group which contains the direct product 
of two (non-abelian) free groups does not admit a solution to the generalized
word problem. 
Recently Rips showed \cite{Rips} - without using Mihailova's result - 
that the double of two free groups of rank $2$ over a normal subgroup 
of index $3$ does not have a generalized word problem solution. 
It was, however, already shown by Gersten \cite{Gersten} 
that such a group must contain $F_2 \times F_2$ and thus Rips's result 
follows from Gersten's and Mihailova's.

The purpose of this paper is to clarify the situation. We will show that
Rips's example is virtually the direct product of two free groups (of finite
rank) and thus Rips's proof yields a new proof for Mihailova's result.
Furthermore we will show that a more general construction of a double
of a group must always contain $F_2 \times F_2$.

\medskip

One might naively think that every finitely presented group which admits
a solution to the word problem but which is neither coherent nor does it 
admit a solution to the generalized word problem must contain 
$F_2 \times F_2$.

A counterexample follows from a construction of 
Rips \cite{Rips2}.
He constructed for every finitely presented group $G$ a short
exact sequence
$$
1 \longrightarrow K \longrightarrow H \longrightarrow G \longrightarrow 1
$$
so that $H$ is a finitely presented word hyperbolic group, and $K$ is - as
a group - finitely generated. (In fact, Rips's construction yields for every
$\lambda$ a finitely presented group $H$ satisfying the small cancellation 
condition $C'(\lambda)$. In particular,
$\lambda$ can be chosen so that $H$ is word hyperbolic.)

It follows that if either $G$ is not coherent or $G$ does not admit a solution to the
generalized word problem, then the same is true for $H$.

Since $H$ is word hyperbolic, it has a solution to the word problem but
cannot contain $F_2 \times F_2$ or even $\Z \times \Z$.

Rips's construction starts with a group $G$ that is already known
to be incoherent and which does not admit a solution to the generalized 
word problem.
Very recently, however, Wise \cite{Wise} 
gave a construction of an incoherent word hyperbolic
group without referring to a group which contains $F_2 \times F_2$.
 
\begin{definition}
For a group $G$ one defines the following properties:
\begin{enumerate} 
\item The group $G$ is coherent if every finitely generated subgroup
is finitely presented. 
\item It is LERF (locally extended residually finite) or subgroup separable 
if for every finitely generated subgroup $H$ of $G$ and every $g \in G, \, 
g\not \in H$ there is a subgroup of finite index in $G$ 
that contains $H$ but does not contain $g$. In particular, a LERF group
is residually finite.
\item  The group $G$ admits a solution to the generalized word problem or
occurrence problem if there is an algorithm that decides for every 
element $g \in G$ and every finitely generated subgroup $H \subset G$ 
whether $g$ is in $H$ or not.
\end{enumerate}
\end{definition}

As for a finitely presented group $G$ residual finiteness implies that
the word problem is solvable, it is easy to see that if $G$ 
is LERF then the generalized word problem for $G$ is solvable. 

Moreover, it follows more or less directly from the definition that if a 
group $G$ has one of these properties then so does every subgroup and 
every finite extension of $G$.

An important obstacle for a group to be coherent, LERF or to admit
a solution to the generalized word problem is given by
\begin{proposition}
The group $F_2 \times F_2$ is neither coherent (\cite{Grunewald1})
nor does it admit a solution to the generalized word problem (Mihailova,
see \cite{Miller} and compare with \cite{Grunewald1}) and is thus not LERF.

Therefore every group which contains $F_2 \times F_2$ has also none
of these properties.
\end{proposition}

\noindent
Let $G$ be a group and $H$ a subgroup.
By the double of $G$ over $H$ we mean the group that we get by
the free product of $G$ with itself amalgamated over $H$.
Here we will show that a large family of group doubles
contain $F_2 \times F_2$ as a subgroup.

\begin{lemma}
\label{lem:free}
  Let $G$ be a group and $H \subgroup G$. Let $L = G*_{H}G$ be the double of G over
  H. Let $\phi_1 : L \arrow G$ be the homomorphism obtained by 
identifying the two
  factors of $L$. Then $K = ker(\phi_1)$ is a free group. If $H$ is of
finite index in $G$ then its rank is $[G:H] - 1$.
\end{lemma}
\begin{proof}
The group L is an amalgam of groups and has a natural action on it's 
Bass-Serre tree $T$
(which is locally finite iff $[G:H] < \infinity$). By definition 
$K$ intersects trivially the two copies of $G$ in $L$, and, since $K$ 
is normal, it intersects trivially all their conjugates, which are all 
the vertex stabilizers. We conclude
that $K$ acts freely and effectively on $T$. 
(Note that $L$ does not necessarily act effectively on $T$.) 
Therefore $K \cong \pi_1(K \verb|\| T,\{pt\})$ is
a free group (see also \cite{Serre2}). Now $\phi_1$ is a surjective 
homomorphism and therefore the index of $K$ in $L$ is
$[L:K] = \#\left|G\right|$. (Note that $G$ might
be infinite.)  The graph $K \verb|\| T$ is, therefore, a
$\#\left|G\right|-$fold cover of the edge of groups corresponding to the 
amalgam
$L$ (in the sense of Bass \cite{Bass}). The only possibility of such a 
covering graph is shown in Figure \ref{fig:index}.
   \begin{figure}
   \begin{center}
   \includegraphics[angle = 0]{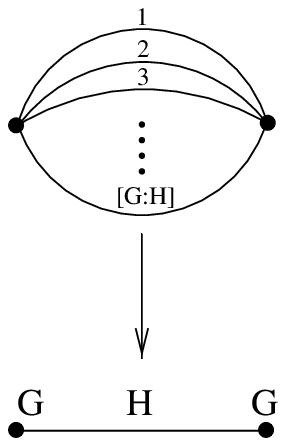}
   \end{center}
   \caption{The covering of graphs corresponding to $K \inclusion L$} \label{fig:index}
   \end{figure}
The rank of $K$ as a free group is now easily verified as $K$ is the fundamental
group of the covering graph in Figure \ref{fig:index}. 
\end{proof}

The idea for the proof of the following theorem was inspired by the notion of 
reducible lattices acting on products of trees, defined by Burger, Mozes
and Zimmer in (\cite{BMZ:Representations,BM:Simple})

\begin{theorem} 
\label{thm:double}
   Let $G$ be a group, $H<G$ be a subgroup with $[G:H] \ge 3$.
   Assume that there exists a subgroup $N \subgroup H$ with the 
   following properties.
   \begin{enumerate}
     \item
     \label{itm:free_subgroup}
        $N$ contains a non-abelian, free subgroup.
     \item
     \label{itm:normal}
        $N$ is normal in $G$.
   \end{enumerate}
   Then the double $L \defeq G*_{H}G$ contains $F_2 \times F_2$ as a subgroup
and is thus neither coherent nor LERF nor does it admit a solution to the generalized
word problem.
\end{theorem}
\begin{proof}
The group $N$ is normal in both copies of $G$ in $L$ and is therefore 
normal in $L$.
We consider the following homomorphism of L into a product:
\begin{equation}
\label{eqn:embeding}
   \phi = (\phi_1,\phi_2) : L  \arrow G \times L/N
\end{equation}
Here $\phi_1$ is induced by identifying the two factors of $L$,
and $\phi_2$ is just the projection to the quotient.
First notice that $ker(\phi_2) = N$ 
and that $N$, as a subgroup of G, injects into G under $\phi_1$. Therefore
$\ker(\phi) = \ker(\phi_1) \cap \ker(\phi_2)$ is trivial. 
So $\phi$ is injective
and it will be enough to prove that $\phi(L)$ contains a copy of 
$F_2 \times F_2$.
Define $K_i \defeq \phi_i(\ker(\phi_{3-i})) \, \forall i \in \{1,2\}$, 
and look at the embedding
$K_1 \times K_2 \inclusion L$. We will show that each one of the $K_i$ contains
a non-abelian free group and thus conclude the theorem. $K_1 = N$ contains a
non-abelian free group by assumption. In order to deal with the other factor, consider the
map $\overline{\phi_1}:G/N*_{H/N}G/N \arrow G/N$ 
(the homomorphism obtained by identifying the two copies of $G/N$).
The fact that $ker(\phi_1) \cap ker(\phi_2) = \left< e \right>$
and the commutativity of the following diagram:

$$\begin{CD}
        L \cong G*_{H}G    @>\phi_{2}>>   
              L/N \cong G/N*_{H/N}G/N              \\
        @V\phi_{1}VV                               @V\overline{\phi_{1}}VV         \\
        G                  @>>>                    G/N                           
\end{CD}
$$

\medskip

\noindent
imply that $K_2 = \ker(\overline{\phi_{1}})$.
Now, by the Lemma \ref{lem:free}, $K_2$ is a free group of rank 
$$[G/N:H/N] - 1 = [G:H] - 1.$$  The assumption that this number is at least 2 
shows that $K_2$ is non-abelian and concludes 
the theorem.
\end{proof}

\begin{corollary}
\label{thm:corollary}
   Let $G$ be a group containing a non-abelian free subgroup.
   Let $H<G$ be a finite index subgroup with $[G:H] \ge 3$. 
   Then the double $L \defeq G*_{H}G$ contains $F_2 \times F_2$ as a subgroup.
\end{corollary}
\begin{proof}
Take $N$ to be a finite index subgroup of $H$ which is normal in $G$, and
all the conditions of Theorem \ref{thm:double} are satisfied.
\end{proof}

\medskip

\noindent
For a property $P$, a group is virtually $P$ if a subgroup
of finite index has this property. 
As another corollary we have:

\begin{corollary}
Let $G$ be a virtually (non-abelian) free group and $H$ a subgroup of 
finite index at least $3$. The double $G \ast_{H} G$ is virtually 
the direct product $F_{r_1} \times F_{r_2}$ of two non-abelian 
free groups of finite rank.

In particular if $G$ is a non-abelian free group then 
$F_2 \times F_2$ and $G \ast_{H} G$ 
have a common subgroup of finite index.
\end{corollary}
\begin{proof}
This follows directly from the proof of Theorem \ref {thm:double}.
In our situation $G \ast_H G$ contains $K_1 \times K_2$ as a subgroup
of finite index where $K_1$ is virtually free and $K_2$ is free. 
\end{proof}

\noindent
{\bf \Large Acknowledgment:}

\smallskip

\noindent
This paper was written while the second and third author were guests
at Columbia University. They would like to thank Columbia and especially
Joan Birman and Hyman Bass for their hospitality. The authors would also 
like to thank Michael Kapovich for bringing to their attention and explaining
to them Rips's construction of incoherent word hyperbolic groups
mentioned in the introduction.

\bigskip

\noindent
{\bf \Large Addresses:}

\smallskip
\noindent
Columbia University, Department of Mathematics, New York, NY 10027

\smallskip
\noindent
Columbia University, Department of Mathematics, New York, NY 10027 and\newline
\noindent
Heinrich-Heine-Universit\"at, Mathematisches Institut,
Abt. M.+M.D., D-40225 D\"usseldorf

\smallskip
\noindent
Institute of Mathematics, The Hebrew University, Jerusalem 91904, Israel

\smallskip
\noindent
Barnard College, Columbia University, 
Department of Mathematics, New York, NY 10027

\bibliographystyle{alpha}
\bibliography{../linklit}

\end{document}